\documentclass{amsart}

\usepackage{graphicx}

\usepackage{amscd}

\usepackage{amssymb}
\usepackage{amsthm}
\newtheorem{theorem}{Theorem}
\newtheorem{lem}[theorem]{Lemma}
\newtheorem{defin}[theorem]{Definition}
\newtheorem{proposition}[theorem]{Proposition}

\newtheorem{conjecture}[theorem]{Conjecture}
\newtheorem{corollary}[theorem]{Corollary}
\numberwithin{equation}{section}

\def\ker#1{{\rm ker}(#1)}
\def\irr#1{{\rm Irr}(#1)}
\def\norm#1{{\rm n}(#1)}
\def\ordo#1{{\rm o}(#1)}
\def\logp#1{{\rm log_p}(#1)}
\usepackage{amsmath}
\begin{document}

\title{COVERINGS OF ABELIAN GROUPS AND VECTOR SPACES}
\author{BAL\'AZS SZEGEDY}
\date{}
\maketitle

\begin{abstract}\noindent
We study the question how many subgroups, cosets or subspaces are needed to cover a finite Abelian group or a vector space
if we have some natural restrictions on the structure of the covering system.
For example we determine, how many cosets we need, if we want to cover all but one
element of an Abelian group. This result is a group theoretical extension of
the theorem of Brouwer, Jamison and Schrijver about the blocking number of an affine space.
We show that these covering problems are closely related to combinatorial problems, including the so called additive basis conjecture,
the three-flow conjecture, and a
conjecture of Alon, Jaeger and Tarsi about nowhere zero vectors.


\end{abstract}

\section{Introduction}
\smallskip
A subgroup covering (coset covering) of the group $G$ is a collection of its subgroups (cosets of its subgroups) whose union is the whole group.
A covering is called irredundant or minimal if none of its memebers can be omitted.
B.H. Neumann observed \cite{NE} that if a group $G$ has a finite irredundant (right) coset covering $H_1x_1~,~ H_2x_2~,~...~,~H_nx_n$
then the index of the intersection of the
subgroups $H_i$ is bounded above by some function of $n$. Let $f(n)$ (resp. $g(n)$) be the maximal possible value
of $|G:\bigcap H_i|$ where $G$ is a group with a coset covering $\{H_ix_i | i=1...n\}$
(resp. subgroup covering $\{H_i | i=1...n\}$).
Obviously we have $f(n)\geq g(n)$. M.J. Tomkinson proved \cite{T} that $f(n)=n!$ and that $g(n)\geq {{1}\over{2}} \cdot 3^{2(n-1)/3}$.
Since no super exponential lower bound has been found for $g(n)$, its order of magnitude is conjectured to be exponential.
Let the functions $f_1(n)$ and $g_1(n)$ be similarly defined as $f(n)$ and $g(n)$ with the additional restriction that the group $G$
is always assumed to be Abelian.
(Note that $f_1(n)\leq f(n)$ and $g_1(n)\leq g(n)$) L. Pyber pointed out (see: \cite{py}) that the order of magnitude of $g_1(n)$ is itself interesting.

We need the following definition.

\begin{defin}
Let $G$ be a fixed finite group. Let $f(G)$ (resp. $g(G)$) denote
the minimal k such that there exists an irredundant covering by
$k$ cosets $\{H_ix_i | i=1...k\}$ (resp. subgroups $\{H_i |
i=1...k\}$) of $G$ where $\bigcap H_i$ is trivial. (Note that the
set of such subgroup coverings may be empty, and in this case we
define $g(G)$ to be infinite)
\end{defin}
Now we have that $g(G)\geq f(G)$. Pyber's problem transforms to find a logarithmic lower bound for $g(A)$ in terms of $|A|$ if
$A$ is an Abelian group.

\begin{conjecture}[Pyber]\label{pybcon}
There exists a fixed constant $c>1$ such that
$g(A)> {\rm log}_c|A|$ for all finite Abelian groups $A$.
\end{conjecture}

Actually we believe that (in contrast with $f(n)=n!$) the growth of $f_1(n)$ is bounded above by some exponential function
and thus

\begin{conjecture}\label{pybcon2}
There exists a fixed constant $c_2>1$ such that
$f(A)> {\rm log}_{c_2}|A|$ for all finite Abelian groups $A$.
\end{conjecture}

We note that the worst known case (even for the function $f(A)$ and thus for $f_1(n)$) is the elementary Abelian $2$-group $A={C_2}^n~(n>1)$ where $f(G)=g(G)=n+1$
(See Corollary \ref{elemi}). It suggests that perhaps $2$ could be the true value for the constant $c$.

We have two results related to Conjecture \ref{pybcon} and Conjecture \ref{pybcon2}.

\begin{theorem}\label{fedthm}
Let $A$ be an Abelian group of order $p_1^{\alpha_1}p_2^{\alpha_2}\dots p_n^{\alpha_n}$ . Then $g(A)\geq f(A)\geq 1+\sum_{i=1}^n \alpha_i$.
\end{theorem}

It means in particular that the inequality of Conjecture \ref{pybcon} holds with $c=p_n$ where $p_n$ is the largest
prime divisor of the order of $A$.

Alon and F\"uredi in \cite{AF} prove the suprising result that if want to cover all but one vertices of an $n$-dimensional cube,
then we need at least $n$ hyperplanes. Actually they prove a more general result.

\smallskip
\noindent{\bf Theorem }(Alon, F\"uredi).
{\it Let $h_1,h_2,\dots,h_n$ be positive integers and let $V$ be the set of all lattice points $(y_1,y_2,\dots,y_n)$ with
$0\leq y_i\leq h_i$. If we want to cover all but one of the points of $V$, then we need at least $h_1+h_2+\dots+h_n$
hyperplanes.}
\smallskip

Our next result is an analogy of the previous one. We determine, how many cosets we need, if we want to cover all but one
element of an Abelian group.
This result yields good lower bound for the size of an irredundant coset covering system if it contains a small coset.

\begin{theorem}\label{mainthm}
Let $A$ be an Abelian group of order $p_1^{\alpha_1}p_2^{\alpha_2}\dots p_n^{\alpha_n}$.
Let $\phi(A)$ denote the minimal number $k$ for which there exists a system of subgroups $H_1$,$H_2$,\dots,$H_k$ and
elements $x_1,x_2,\dots,x_k$ such that $G\setminus\{1\}=\bigcup H_ix_i$.
Then $\phi(A)=\sum_1^n \alpha_i (p_i-1)$.
\end{theorem}
\smallskip
\begin{corollary}\label{mainthmc}
Let $A$ be an Abelian group and
let $\{H_ix_i | i=1...k\}$ be an irredundant coset covering of $A$. Then for all $i$
\begin{equation*}
k\geq 1+{\rm log}_2|G:H_i|
\end{equation*}
\end{corollary}
Note that Theorem \ref{mainthm} solves the special case of conjecture \ref{pybcon2} when one of the cosets
has size 1. In this case both conjectures hold with constant 2. Corollary \ref{mainthmc} shows that if the covering
system contains a "small" subgroup of size less than $|A|^p$ for some $p<1$ then both conjectures hold with
constant $c=c_2=2/(1-p)$.

The proof of Theorem \ref{mainthm}. uses character theory and some
Galois theory. It is also worth mentioning that Theorem
\ref{mainthm} implies that the blocking number of an affine space
(i.e. the size of the smallest subset which intersect all
hyperplanes) over the prime field GF($p$) is $1+n(p-1)$ which was
proved (for arbitrary finite fields) by Brouwer, Schrijver and
Jamison  \cite{BS}, \cite{jamison} using polynomial method.

\medskip
\smallskip
From the combinatorial point of view, the most important special
case of the previously described covering problems is when the
group $A$ is an elementary Abelian group $(C_p)^n$, and the
covering system consists of hyperplanes (or affine hyperplanes).
More generally, we can speak about hyperplane coverings of vector
spaces over arbitrary finite fields. Many questions about graph
colorings, nowhere zero flows or nowhere zero vectors can be
translated to questions about special hyperplane coverings.
However not much is known about such coverings. In Chapter 5. we
present a character theoretic approach to hyperplane coverings in
vector spaces over prime fields. The space of n-dimensional row
vectors admits a natural scalar product. We prove the following:

\begin{theorem}
Let $p$ be an odd prime and let $A={\rm GF}(p)^n$.
The space $A$ is covered by the hyperplanes
${{\mathbf x}_1}^{\bot},{{\mathbf x}_2}^{\bot},\dots,{{\mathbf x}_k}^{\bot}$
if and only if for all vectors ${\mathbf v}\in A$ the number of 0-1 combinations of the vectors
${{\mathbf x}_1},{{\mathbf x}_2},\dots,{{\mathbf x}_k}$ resulting ${\mathbf v}$ is even.
\end{theorem}

\medskip
\noindent{\bf Conjecture} (Alon, Jaeger, Tarsi).
{\it Let $F$ be a finite field with $q>3$ elements and let $M$ be a nonsingular $n$ by $n$ matrix over $F$.
Then there exists a nowhere zero (column) vector $x$ (i.e. each component of $x$ is non zero) such that the vector $Mx$
is also a nowhere zero vector.}
\smallskip

With an elegant application of the polynomial method of Alon, Nathanson and Ruzsa (see: \cite{ANR})
Alon and Tarsi prove \cite{AT} that the latter conjecture holds if $F$
is a proper extension of some prime field GF($p$).
Actually they prove more. For example, from their results follows that if ${\mathbf v}$ is an arbitrary (column) vector
and $M$ is a nonsingular matrix over $F$ then there exists a nowhere zero vector ${\mathbf x}$ such that
$M{\mathbf x}-{\mathbf v}$ is a nowhere zero vector. It is reasonable to believe that the same statement
holds over GF($p$) where $p$ is prime number and is bigger than $3$. This conjecture will be called the
{\bf choosability} version of the Alon-Jaeger-Tarsi conjecture.

\begin{proposition}\label{pimpa}
A positive answer of conjecture \ref{pybcon}
implies the Alon-Jaeger-Tarsi conjecture for $F$=GF($p$) where $p\geq c^2$.
\end{proposition}
 In Chapter 7.
we discuss minimal hyperplane coverings.

\begin{defin}
Let $V$ be an $n$-dimensional vectorspace over the finite field $GF(q)$. Let $h_q(n)$ ($l_q(n)$) denote
the minimal number $k$ such that there is collections of $k$ hyperplanes (affine hyperplanes)
of $V$ which forms a minimal covering system and the intersection of these hyperplanes
(the hyperplanes corresponding to these affine hyperplanes) is trivial.
\end{defin}

Obviously we have $h_q(n)\geq l_q(n)>n$.
Let $$1+\varepsilon_q=\inf_n{l_q(n)/n}.$$
We conjecture that $\varepsilon_q>0$ if $q$ is an arbitrary prime power bigger than 2.
This conjecture can be formulated in the following nice, self contained way.

\begin{conjecture}\label{thecon}
Assume that for some prime power $q>2$ the $GF(q)$ vectorspace $V$ is
covered irredundantly by $k$ affine hyperplanes $H_1+v_1,H_2+v_2,\dots H_k+v_k$. Then
the codimension of the intersection $\bigcap_i H_i$ is at most $k/(1+\varepsilon_q)$ for some fixed
positive constant $\varepsilon_q$ which depends only on $q$.
\end{conjecture}

Using a result of Alon and Tarsi about nowhere zero points \cite{AT}, we prove the following.

\begin{theorem}
If $q$ is not a prime number then $\varepsilon_q\geq {1\over 2}$.
\end{theorem}

The $p=3$ case of Conjecture \ref{thecon} is especially interesting because it is strongly related to the next
two conjectures.

\smallskip
\noindent{\bf Weak three flow conjecture}
{\it There exits a fixed natural number $k$ such that if a graph $G$ is at least $k$-connected
then it admits a nowhere zero $3$-flow.}
\smallskip

It is well known that the next conjecture (for $p=3$) would imply the weak 3-flow conjecture.

\smallskip
\noindent{\bf Additive basis conjecture}\quad (Jaeger, Linial,
Payan, Tarsi). {\it For every prime $p$ there is a constant $c(p)$
depending only on $p$ such that if $B_1,B_2,\dots ,B_{c(p)}$ are
bases of the $GF(p)$ vectorspace $V$ then all elements of $V$ can
be written as a zero-one linear combination of the elements of the
union (as multisets) of the previous bases.}
\smallskip

We show that $\varepsilon_3>0$ is equivalent with the additive basis conjecture for $p=3$.
For a prime number $p>3$ we show that $\varepsilon_p>1$ implies the choosability version of the
Alon-Jaeger-Tarsi conjecture and that the latter one implies $\varepsilon_p\geq 0.5$.
Note that Conjecture \ref{pybcon2} implies that $\varepsilon_p\geq {\log}_2p-1$.

\bigskip

\section{Notation and basics}
\smallskip

Let $A$ be a finite Abelian group. A linear character of $A$ is a
homomorphism from $A$ to $\Bbb{C}^*$. The linear characters of $A$
are forming a group under the point wise multiplication (which is
isomorphic to $A$) and they are forming a basis in the vector
space of all function $f:A\rightarrow \Bbb{C}$. The trivial
character (which maps all elements of $A$ to $1$) will be denoted
by $1_A$. The kernel of a linear character $\chi$ is the set of
those group elements $g\in A$ for which $\chi(g)=1$. We denote the
kernel of $\chi$ by ${\rm ker}(\chi)$. It is easy to see that the
subgroup $H\leq A$ is the kernel of some liner character $\chi$ if
and only if $A/H$ is cyclic.

The group algebra  $\Bbb{C}A$ consists of the formal linear
combinations of the group elements. The fact that some Abelian
groups are imagined as additive structures can cause some
confusion because the concept of group algebra suggests that the
group operation is the multiplication. For example we will work in
the group algebra of the additive group of a finite vector space
$V$. In this structure all vectors from $V$ are linearly
independent and the group algebra $\Bbb{C}V^+$ consists of the
formal $\Bbb{C}$-linear combinations of the elements of $V$. The
product of two vectors ${\mathbf v}_1$ and ${\mathbf v}_2$ is the
vector ${\mathbf z}={\mathbf v}_1+{\mathbf v}_2$ with coefficient
$1$. If we add together ${\mathbf v}_1$ and ${\mathbf v}_2$ in the
group algebra then it has nothing to do with the element ${\mathbf
z}$. Another source of confusion is that the identity element of
the group algebra is the zero vector with coefficient one. The
identity element of the group algebra is always denoted by $1$.
For a good reference about characters and group algebras see
\cite{i}.
\smallskip

Let $V$ be an $n$ dimensional vector space.
A Hyperplane of $V$ is a subspace of co-dimension 1.
We say that the hyperplanes $H_1,H_2,\dots,H_k$ are independent
(or the set $\{H_1,H_2,\dots,H_k\}$ is independent)
if the co-dimension of their intersection is $k$.
If $V$ is represented as the space of row vectors of length $n$ then
there is a natural scalar product on $V$ defined by $({\mathbf x},{\mathbf y})=\sum_{i=1}^n x_iy_i$.
The vectors ${\mathbf x}_1,{\mathbf x}_2,\dots,{\mathbf x}_k~\in V$ are linearly independent if and only if the hyperplanes
${{\mathbf x}_1}^{\bot},{{\mathbf x}_2}^{\bot},\dots,{{\mathbf x}_k}^{\bot}$ are linearly independent.
If $V$ is a row space, then the usual basis will be always denoted by
${\mathbf b}_1,{\mathbf b}_2,\dots$, their orthogonal spaces will be denoted
by $B_1,B_2,\dots$ and we call them basis hyperplanes.
An affine hyperplane is the set $H+{\mathbf v}$ where $H$ is a hyperplane and ${\mathbf v}$ is a vector.
If $A=H+{\mathbf v}$ is an affine hyperplane, we say that $H$ is the hyperplane corresponding to $A$.
A collection of affine hyperplanes $A_1,A_2,\dots,A_k$ is called independent if the corresponding hyperplanes
are independent.
\bigskip
\section{Proof of Theorem 4.}

Let $\Omega=\{H_1x_1,H_2x_2,\dots,H_kx_k\}$ be a coset system of
the Abelian group $A$. We say that $\bigcap_{i=1}^k H_i$ is the
{\bf subgroup intersection} of $\Omega$. If $S$ is a subset of
$A$, we say that $S$ is {\bf covered} by $\Omega$ if it is
contained in the union of the elements of $\Omega$. Let $M$ be a
subgroup of $A$. We denote by $\Omega/M$ the coset system in $A/M$
consisting of the images of the cosets
$H_1x_1,H_2x_2,\dots,H_kx_k$ under the homomorphism $A\rightarrow
A/M$. By abusing the notation, we denote by $\Omega\cap M$ the
system consisting of the cosets $$H_1x_1\cap M,~H_2x_2\cap
M,~\dots,~H_kx_k\cap M.$$

\begin{proof}[Proof of theorem 4.]
For a natural number $n$ with prime decomposition
$p_1^{\alpha_1}p_2^{\alpha_2}\dots p_m^{\alpha_m}$ let
$\lambda(n)=\alpha_1+\alpha_2+\dots+\alpha_m$. We prove Theorem
\ref{fedthm}. by induction on the order of the Abelian group
$|A|$. During the proof we will frequently use the fact that the
coset structure of $A$ is translation invariant. Let
$\Omega=\{H_1x_1,H_2x_2,\dots,H_kx_k\}$ be a minimal coset
covering system of $A$ with trivial subgroup intersection. We have
to show that $k\geq 1+\lambda(|A|)$.

Let $K$ be a maximal subgroup of $A$ containing $H_1$. Note the
$K$ has prime index in $A$ and so $\lambda(|A|)=\lambda(|K|)+1$.
Using the fact that any translation of the system $\Omega$ is
again a minimal coset covering system with trivial subgroup
intersection, we can assume that $x_1\notin K$. This means that
$H_1x_1$ is disjoint from $K$ and that $H_2x_2,\dots,H_kx_k$
covers $K$. Let $\Omega_1\subseteq
\{H_2x_2,\dots,H_kx_k\}\subset\Omega$ be minimal with the property
that it covers $K$. There are two possibilities.

The first one is that the subgroup intersection of $\Omega_1$ is
trivial. In this case the subgroup intersection of $\Omega_1\cap
K$ is also trivial and then by induction we can deduce that
$$k\geq 1+f(K)\geq 2+\lambda(|K|)=1+\lambda(|A|)$$ which finishes
the proof.

The second one is that the subgroup intersection of $\Omega_1$ is
not trivial. let $M_1$ denote the subgroup intersection of
$\Omega_1$. Since the factor group $K/(K\cap M_1)$ is covered
minimally by $\Omega_1/M_1$ with trivial subgroup intersection, we
have by induction that $|\Omega_1/M_1|=|\Omega_1|\geq
1+\lambda(K/(K\cap M_1))$. Let $y_1$ be an element of $A$ which is
not covered by the cosets in $\Omega_1$. Clearly the whole coset
$M_1y_1$ does not intersect any coset from $\Omega_1$. Let
$\Omega_2\subseteq\Omega$ be a minimal covering system for
$M_1y_1$ and let $M_2$ be the subgroup intersection of
$\Omega_1\cup\Omega_2$. Using translation invariance we have that
$\Omega_2\geq 1+\lambda(M_1/M_2)$.

\smallskip
Now we define a process. Assume that the $\Omega_i$, $M_i$ and
$y_i$ is already constructed for $1\leq i\leq t$ and the subgroup
intersection $M_i$ is still not trivial. Let
$\Omega_{t+1}\subseteq\Omega$ be a minimal covering for $M_ty_t$.
Let $M_{t+1}$ denote the subgroup intersection of the system
$\bigcup_{i=1}^{t+1}\Omega_i$. If $M_{t+1}$ is not trivial then
let $y_{i+1}$ be an element which is not covered by the system
$\bigcup_{i=1}^{t+1}\Omega_i$.

Using the induction hypothesis and translation invariance we get
that $|\Omega_{t+1}|\geq 1+\lambda(M_t/M_{t+1})$. Assume that
$M_r$ is trivial, and thus $r$ is the length of the previous
process. Now we have that
$$|\Omega|\geq\sum_{i=1}^r |\Omega_i|\geq r+\lambda (|K|)\geq
1+\lambda(|A|).$$
\end{proof}

Using Theorem 1. we obtain precise result for $(C_2)^n$.

\begin{corollary}\label{elemi}
$f((C_2)^n)=n+1$
\end{corollary}

\begin{proof}
Theorem 1. implies that $f((C_2)^n)\geq n+1$. Let $H_i$~($1\leq i\leq n$) be the subgroup consisting of
all elements with 0 at the i-th place. The group $(C_2)^n$ is the union of the groups $H_i$ and the
element $(1,1,\dots,1)$.
\end{proof}

\section{Proof of Theorem 5.}
\smallskip
For a natural number $n$ with prime decomposition
$n=\prod_{i=1}^{s} {p_i}^{\alpha_i}$, let
$\tau(n)=\sum_{i=1}^{s}{\alpha_i}(p_i-1)$. Let $\phi(A)$ denote
the smallest number $k$ for which there is a collection of cosets
$H_1x_1,H_2x_2,\dots,H_kx_k$ in the Abelian group $A$ such that
$$A\setminus \{1\}=\cup_{i=1}^k H_ix_i.$$

\begin{lem}\label{lem0}
        $\phi(A)\leq\tau(|A|)$
\end{lem}

\begin{proof}
We go by induction on $|A|$. Let $B<A$ be a subgroup of index $p_1$.
The group $A$ is a disjoint union of $p_1$ cosets of $B$.
Using the statement for $B$ we obtain the result for $A$.
\end{proof}

\begin{lem}\label{lem1}
       Let $B$ and $C$ be two Abelian groups with $(|B|,|C|)=1$.
Then $\phi(B\times C)\geq\phi (B)+\phi (C)$.
\end{lem}

\begin{proof}
If two groups have coprime order then a subgroup of their direct product
is a direct product of their subgroups.
It follows that for $H\leq B\times C$ and $g\in B\times C$ there are
subgroups $H_1\leq B$ , $H_2\leq C$ and elements $g_1\in B$ , $g_2\in C$
such that $Hg=\{(h_1g_1,h_2g_2)|h_1\in H_1 , h_2\in H_2\}$.
Assume that $(B\times C)\setminus\{1\} = \bigcup_1^k{K_ig_i}$
where $K_i<B\times C$, $g_i\notin K_i$ and $k=\phi (B\times C)$.
If $K_ig_i$ intersects $(B,1)\leq B\times
C$ then it does not intersect $(1,C)$ otherwise
the 1 would be an element of $K_ig_i$. (The analogous statement holds if
$K_ig_i$ intersects $(1,C)$.)
This implies that $k\geq \phi (B)+\phi (C)$.

\end{proof}

\begin{lem}\label{lem3}
       If $H<G$ and $g\notin H$ then there exists a subgroup
$K$ of $G$ such that $H\leq K$, $g\notin K$, $G/K$ is cyclic.
\end{lem}

\begin{proof}
Let $K$ be maximal with the property $H\leq K <G$,
$g\notin K$. In the factor group $G/K$ every nontrivial
subgroup $K_2$ contains $Kg$ otherwise the preimage
of $K_2$ under the homomorphism $G\rightarrow G/K$ would be
bigger than $K$ and would not contain $g$.
It follows that $G/K$ can't be a direct product
of two proper subgroups because one of them would not contain $Kg$.
Using the structure theorem of finite Abelian groups we obtain
that $G/K$ must be cyclic of prime power order.
\end{proof}

\begin{lem}\label{lem4}
If $P$ is an Abelian group of order $p^{\alpha}$ for some
prime $p$ and integer $\alpha$, then $\phi(P)\geq \alpha (p-1)$.
\end{lem}

\begin{proof}
Let $k=\phi(P)$ and $P\setminus \{1\}=\bigcup_1^k{H_i}g_i$ (where
$g_i\notin H_i$). Using Lemma \ref{lem3}. we obtain that there are subgroups $K_i$ $(1\leq i \leq k)$
with $H_i\leq K_i$, $g_i\notin K_i$ and $P/K_i$ is cyclic for all $1\leq i\leq k$.
Now we have $P\setminus\{1\}=\bigcup_1^k{K_i}g_i$ and
for each $K_i$ there exists a linear character $\chi_i$ of $P$
such that $\ker{\chi_i}=K_i$.
Clearly the product $\prod_1^k{(\chi_i-(\chi_i(g_i))1_P)}$ takes
the value zero
on every element $1\neq g\in P$ but it is nonzero on the element $1$.
From this we obtain the following equality
$$\prod_1^k{(\chi_i-\chi_i(g_i)1_P)}=
(\prod_1^k{(1-\chi(g_i))}/|P|)(\sum_{\chi \in \irr{P}}\chi)$$
The linear characters of $P$ are forming a basis of the vector space
of $P\rightarrow \Bbb{C}$ functions, and thus after
expanding both sides of the above equation, the coefficients
of the characters must coincide.
On the left hand side each coefficient is a sum of roots of unities thus
they are algebraic integers.
On the right hand side every character has coefficient
$\prod_1^k{(1-\chi(g_i))}/|P|$, and thus this number is an algebraic integer.
The $|P|$-th cyclotomic field $F$ is a normal extension
of $\Bbb{Q}$, and the degree of the field extension
$F/\Bbb{Q}$ is $p^{\alpha-1}(p-1)$.
Using the fact that the Galois norm of an algebraic
integer is an integer we deduce that $\norm{|P|}=$
$p^{\alpha p^{\alpha-1}(p-1)}$ divides $\prod_1^k{\norm{1-\chi(g_i)}}$
where $\norm{x}$ denotes the Galois norm of $x$ in the
field extension $F/\Bbb{Q}$.
An easy calculation shows that
$\norm{1-\chi(g_i)}=p^{p^{\alpha-\logp{\ordo{\chi(g_i)}}}}\leq$
$p^{p^{\alpha-1}}$
where $\ordo{\chi(g_i)}$ denotes the multiplicative order of $\chi(g_i)$.
The last inequality completes the proof.
\end{proof}

\begin{proof}[Proof of Theorem 5.]
  According to Lemma \ref{lem0}, it is enough to prove that $$\phi(A){\geq} \tau(|A|).$$
  We go by induction on $|A|$. If $|A|$ is a prime power then
Lemma \ref{lem4} yields the result.
If $A$ is not a prime power then $A=B\times C$ where $(|B|,|C|)=1$ and
using the statement for $B$ and $C$, Lemma \ref{lem1}
completes the proof.
\end{proof}


\begin{proof}[proof of Corollary \ref{mainthmc}]
Let $g$ be an element of $H_ix_i$ which is not covered by $H_jx_j$ for all $j\neq i$.
Lemma \ref{lem0} shows that there is a coset system $\Omega$ consisting of $\tau(|H_i|)$ cosets whose
union is $H_ix_i\setminus\{g\}$.
The union of the system
$\Omega\cup\{H_jx_j|j\neq i, 1\leq j\leq k\}$ is $A\setminus \{g\}$, so translating
it with $g^{-1}$ we can apply Theorem \ref{mainthm}. We obtain that $k-1+\tau(|H_i|)\geq \tau(|A|)$ and thus
$k\geq 1+\tau(|A:H_i|)$.
It means in particular that $k\geq 1+{\rm log}_2|G:H_i|$.
\end{proof}
\medskip
\section{Hyperplane coverings and characters}
\smallskip

Now we describe our character theoretical approach to hyperplane
covering problems. Let $p$ be a fixed prime number, let
$\omega=e^{2\pi i/p}$ and let $A=(C_p)^n$. We regard $A$ as the
$n$-dimensional row vector space over GF($p$).

\begin{lem}\label{covhyp}
the space $A$ is covered by the hyperplanes ${{\mathbf
x}_1}^{\bot},{{\mathbf x}_2}^{\bot},\dots,{{\mathbf x}_k}^{\bot}$
if and only if the equation
\begin{equation*}
({\mathbf x}_1-{\mathbf 1})({\mathbf x}_2-{\mathbf
1})\dots({\mathbf x}_k-{\mathbf 1})=0.
\end{equation*}
is satisfied in the group algebra ${\Bbb C}[A]$ where ${\mathbf
1}$ denote the identity element of $A$ (which is actually the zero
vector, if we think of $A$ as a vector space).
\end{lem}

Note that substraction in the previous lemma is the group algebra
substraction and not the vector substraction.

\begin{proof}
The function
\begin{equation*}
f:(x_1,x_2,\dots,x_n)\rightarrow((y_1,y_2,\dots,y_n)\rightarrow \omega^{x_1y_1+x_2y_2+\dots+x_ny_n})
\end{equation*}
gives an isomorphism between $A$ and its character group $A^*$. Moreover $f$ can be uniquely extended to an algebra
isomorphism between the group algebra ${\Bbb C}[A]$ and the character algebra ${\Bbb C}[A^*]$.
Note that the character algebra is just the algebra of all functions $A\rightarrow {\Bbb C}$
with the point wise multiplication.
Clearly we have that ${\mathbf x}^\bot = {\rm ker}(f({\mathbf x}))$ for all row vectors ${\mathbf x}\in A$.
It follows that
the space $A$ is covered by the hyperplanes
${{\mathbf x}_1}^{\bot},{{\mathbf x}_2}^{\bot},\dots,{{\mathbf x}_k}^{\bot}$
if and only if
\begin{equation*}
(f({\mathbf x}_1)-1_A)(f({\mathbf x}_2)-1_A)\dots(f({\mathbf x}_k)-1_A)=0.
\end{equation*}
Applying $f^{-1}$ to both side of the previous equation we obtain
the statement of the lemma .
\end{proof}

The previous lemma gives a characterization of covering systems in
terms of orthogonal vectors. Our following theorem gives a group
algebra free characterization of coverings in terms of orthogonal
vectors if $p$ is an odd prime.

\begin{theorem}\label{chariz}
Let $p$ be an odd prime and let $A=(C_p)^n$.
The space $A$ is covered by the hyperplanes
${{\mathbf x}_1}^{\bot},{{\mathbf x}_2}^{\bot},\dots,{{\mathbf x}_k}^{\bot}$
if and only if for all vectors ${\mathbf v}\in A$ the number of 0-1 combinations of the vectors
${{\mathbf x}_1},{{\mathbf x}_2},\dots,{{\mathbf x}_k}$ resulting ${\mathbf v}$ is even.
\end{theorem}

\begin{proof}
Let $F$ be the algebraic closure of the field with two elements. Since $p$ is odd, $F$ contains a $p$-th root
of unity $\omega$ and and thus we can repeat everything
what we did over $\Bbb{C}$. We obtain that the space $A$ is covered by the hyperplanes
${{\mathbf x}_1}^{\bot},{{\mathbf x}_2}^{\bot},\dots,{{\mathbf x}_k}^{\bot}$
if and only if the equation
\begin{equation*}
({\mathbf x}_1-{\mathbf 1})({\mathbf x}_2-{\mathbf 1})\dots({\mathbf x}_k-{\mathbf 1})=0.
\end{equation*}
holds in the group algebra $F[A]$.
Since $F$ has characteristic 2 we don't have to care about the signs in the previous formula.
The rest of the proof is straightforward by expanding the formula.
\end{proof}

\section{on the Alon-Jaeger-Tarsi conjecture}
\smallskip
The following lemma shows the relationship between hyperplane coverings and the Alon-Jaeger-Tarsi conjecture.
\begin{lem}\label{ketfug}
Let $p$ be a fixed prime number and let $n$ be a fixed natural number. The following statements are equivalent.
\begin{enumerate}
\item The $n$ dimensional vector space over GF($p$) can't be covered
by the union of two independent sets of hyperspaces.
\item If $M$ is a non singular $n$ by $n$ matrix over GF($p$) then there exists a nowhere zero vector ${\mathbf x}$ such that
$M{\mathbf x}$ is also a nowhere zero vector.
\end{enumerate}
\end{lem}

\begin{proof}
(1)$\Rightarrow$ (2)
 Let ${\mathbf x}_1,{\mathbf x}_2,\dots,{\mathbf x}_n$ denote the rows of $M$, and let
$H_i={{\mathbf x}_i}^{\bot}$ for $1\leq i \leq n$.
Since $M$ is non singular we have that the subspaces
$H_1,H_2,\dots,H_n$ are independent.
Let $S_i$ be the hyperspace consisting of the row vectors with a zero at the $i$-th component.
It follows from (1) that there exists a vector ${\mathbf y}$ which is not contained in the union of the
spaces $H_i$ , $S_i$~ $(1\leq i\leq n)$. Clearly ${\mathbf y}$ is a nowhere zero vector such that $M{\mathbf y}^T$ is also
a nowhere zero vector.

(2)$\Rightarrow$ (1)

Assume that $V$ is an $n$ dimensional vector space covered by the independent hyperspace sets $\Omega_1$ and $\Omega_2$.
We can assume that both $\Omega_1$ and $\Omega_2$ are maximal independent sets.
It is easy to see that we can represent $V$ as a row space such that the hyperspaces in $\Omega_1$ are exactly
the spaces formed by all vectors with a zero at a fixed component. Let ${\mathbf x_1},{\mathbf x_2},\dots,{\mathbf x_n}$
be a system of non zero vectors whose orthogonal spaces are exactly the hyperspaces in $\Omega_2$.
It is clear that the vectors ${\mathbf x_i}~ (1\leq i\leq n)$ are linearly independent.
Let $M$ be a matrix such that its row vectors are ${\mathbf x_i}$~$1\leq i\leq n$.
Now $M$ contradicts the assumption of (2).
\end{proof}

\begin{proof}[Proof of proposition \ref{pimpa}]
Using Lemma \ref{ketfug} it is enough to show that if $V$
is covered by two sets of independent hyperspaces then $p<c^2$.
Let $\Omega$ be the union of two independent hyperplane sets.
Clearly $\Omega$ contains an independent set $\Delta$ of cardinality $k\geq |\Omega|/2$. Let $W$ denote
the intersection of the hyperspaces in $\Omega$. Now, the factor space $V/W$ is covered irredundantly by the elements
of $\Omega$ and the intersection of this covering system is trivial (in $V/W$). We also have that $d={\rm dim}V/W\geq k$.
It follows that $|\Omega|\leq 2d$.
If conjecture \ref{pybcon} holds then ${\rm log}_c(p^d)< 2d$ which means $p<c^2$.
\end{proof}

\begin{defin}
Let $M$ be an $n$ by $n$ matrix. We say that $M$ is an {\bf AJT matrix} if
there is a nowhere zero (column) vector ${\mathbf x}$ such that $M{\mathbf x}$ is also a nowhere zero vector.
\end{defin}

Note that $M$ is not an AJT if and only if the orthogonal spaces of the
rows of $M$ cover all nowhere zero vectors.

\begin{lem}\label{fedequ}
Let $M$ be an $n$ by $n$ matrix over the field GF($p$) and let $\{{\mathbf x_1},{\mathbf x_2},\dots,{\mathbf x_n}\}$
be the rows of $M$. Moreover let ${\mathbf b_i}$ be the $i$-th row
of the $n$ by $n$ identity matrix. Then $M$ is an AJT if and only if
\begin{equation*}
({\mathbf b}_1-{\mathbf 1})({\mathbf b}_2-{\mathbf 1})\dots({\mathbf b}_n-{\mathbf 1})
({\mathbf x}_1-{\mathbf 1})({\mathbf x}_2-{\mathbf 1})\dots({\mathbf x}_m-{\mathbf 1})\neq 0.
\end{equation*}
in the group algebra ${\Bbb C}[V^+]$ where $V$ denote the space of $n$ dimensional row vectors.
If $p$ is odd and $F$ is the algebraic closure of the field with two elements then the same statements holds if
we replace ${\Bbb C}$ by $F$.
\end{lem}

\begin{proof}
The proof is straightforward from Lemma \ref{covhyp} and Theorem \ref{chariz}.
\end{proof}

\begin{defin}
Let $B$ be a subset of the $n$ dimensional GF($p$) space $V$. Let $C(B)$ denote the set of all vectors ${\mathbf v}$
for which the number of zero-one combinations of the elements from $B$ resulting ${\mathbf v}$ is odd.
In particular if $B$ is a linearly independent set then $C(B)$ is the set of all zero-one combinations
of the elements from $B$. We say that $C(B)$ is the {\bf cube} determined by the set $B$.
Let $A_1,A_2,\dots A_n$ be two element subsets of GF($p$). We say that the vector set $\{(a_1,a_2,\dots,a_n)|a_i\in A_i\}$ is
a {\bf combinatorial cube} in the $n$-dimensional row space.
\end{defin}

Using our character theoretic approach we obtain the following characterization of AJT-s.

\begin{theorem}
Let $M$ be an $n$ by $n$ matrix over the field GF($p$) where
$p>2$. Let $X$ be the set formed by the rows of $M$ and let $B$ be
the ordinary basis of the $n$ dimensional row-space. $M$ is an AJT
if and only if the set $C(X)\cap (C(B)+{\mathbf v})$ has odd
number of points for some vector ${\mathbf v}$.
\end{theorem}

\begin{proof}
Let $F$ be the algebraic closure of the field with two elements.
Recall that the elements of the group algebra $F[V^+]$ are formal $F$-linear combinations of the group elements.

Using Lemma \ref{fedequ} and that $1=-1$ in characteristic 2 we get that M is an AJT if and only if

\begin{equation*}
({\mathbf b}_1+{\mathbf 1})({\mathbf b}_2+{\mathbf 1})\dots({\mathbf b}_n+{\mathbf 1})
({\mathbf x}_1+{\mathbf 1})({\mathbf x}_2+{\mathbf 1})\dots({\mathbf x}_m+{\mathbf 1})=
\end{equation*}

\begin{equation*}
\sum_{S_1\subseteq \{1,2,\dots,n\}}~\sum_{S_2\subseteq \{1,2,\dots,m\}}~
\prod_{i\in S_1}{\mathbf b}_i~\prod_{i\in S_2}
{\mathbf x}_i
\end{equation*}

\noindent is not zero in the group algebra $F[V^+]$.
Let ${\mathbf y}$ be a fixed vector in $V$.
To determine the coefficient of ${\mathbf y}$ in the previous product we have to compute the number of the solutions
of the following equation in $F[V^+]$ where $S_1\subseteq \{1,2,\dots,n\}$, ~$S_2\subseteq \{1,2,\dots,m\}$.
\begin{equation*}
\prod_{i\in S_1}{\mathbf b}_i~\prod_{i\in S_2}{\mathbf x}_i={\mathbf y}
\end{equation*}

Since it does not contain any addition, it can be translated into the following equation in $V$.
\begin{equation*}
\sum_{i\in S_1}{\mathbf b}_i~+\sum_{i\in S_2}{\mathbf x}_i={\mathbf y}
\end{equation*}

The number of the solutions of the previous equation is clearly
$$|C(X)\cap(-C(B)+{\mathbf y})|=|C(X)\cap(C(B)-(1,1,...,1)+{\mathbf
y}|$$ and the parity of this number gives the coefficient of
${\mathbf y}$. It follows that $M$ is an AJT if and only if there
is a vector ${\mathbf v}$ for which $|C(X)\cap C(B)+{\mathbf v}|$
is an odd number.
\end{proof}

As a consequence of the previous lemma we obtain the following.

\begin{corollary}
Let $M$ be an $n$ by $n$ matrix, and let $X$ be the set formed by the rows of $M$.
Then $M$ is an AJT If and only if there is a combinatorial cube which has odd intersection with $C(X)$.
\end{corollary}

\begin{proof}
It is clear that if $M$ is an AJT and $N$ is obtained from $M$ by
multiplying the rows by non zero scalars then $N$ is an AJT too.
Applying the previous theorem to all possible such $N$ the proof
is straightforward.
\end{proof}

\medskip
\medskip

\bigskip
\section{minimal hyperplane coverings}

\begin{lem}\label{brumm}
Let $q$ be a prime power which is not a prime. Let $V$ be an $n$ dimensional vector space over GF($q$)
and let $B_1$ and $B_2$ be two bases of $V$. Then each vector ${\mathbf v}\in V$ can be written as a nowhere zero
linear combination (i.e. neither coefficient is zero) of $B_1\cup B_2$.
\end{lem}

\begin{proof}
We write each vector as a row vector in the basis $B_1$. Let $M$ be a matrix whose rows are the vectors from $B_2$.
According to the results of Alon and Tars in \cite{AT} there is a nowhere zero (row) vector ${\mathbf x}$
such that ${\mathbf v}-{\mathbf x}M$ is a nowhere zero vector ${\mathbf y}$. It means that
${\mathbf v}={\mathbf x}M+{\mathbf y}$ which yields the required linear combination.
\end{proof}

\begin{lem}\label{mats}
Let $M$ be a matroid on the set $E$. If $|E|\geq r(E)k$ for a natural number $k$ then there is a subset $X\subseteq E$
such that $X$ as
a matroid has $k$ disjoint bases.
\end{lem}

\begin{proof}
Let $X$ be a minimal subset of $E$ with the property $|X|\geq r(X)k$. According to Edmond's matroid packing theorem,
the maximal number of pairwise disjoint bases in $X$ equals
$${\rm min}\left\{\left\lfloor{{|X|-|Y|}\over r(X)-r(Y)}\right\rfloor~:~Y\subseteq X~,r(Y)<r(X)\right\}.$$
The minimality of $X$ implies that for an arbitrary subset $Y\subset X$ with $r(Y)<r(X)$ we have that
$|Y|< r(Y)k$. It follows that $|Y|-r(Y)k<|X|-r(X)k$ and so $(|X|-|Y|)/(r(X)-r(Y))>k$.
\end{proof}

\begin{theorem}
Let $q$ be a prime power which is not a prime. Let $V$ be a vector space over GF($q$) which is covered
irredundantly by $k$ affine hyperplanes $H_i+{\mathbf v}_i ~ (1\leq i\leq k)$. Then the co-dimension of the intersection
of the hyperplanes $H_i~(1\leq i\leq k)$ is less than ${2\over 3}k$.
\end{theorem}

\begin{proof}
It is easy to see that the space $V/\bigcap_{1\leq i\leq n} H_i$ is covered
irredundantly by the images of $H_i+{\mathbf v}_i$
so we can assume that $\bigcap H_i$ is trivial. We go by contradiction. assume that ${\rm dim}(V)=n\geq {2\over 3}k$.
Without loss of generality we can assume that $H_1,H_2,\dots H_n$ are independent hyperplanes, and ${\mathbf v}_i$ is
the zero vector for $1\leq i\leq n$.
We can choose a basis $B=\{{\mathbf b}_i|1\leq i\leq n\}$ such that the previous hyperplanes
are exactly the orthogonal spaces of the basis elements.
Let $W=\bigcap_{i>n} H_i$. From our assumption it follows that ${\rm dim}(V/W)\leq k-n\leq {1\over 2}n$.
Let ${\mathbf p}_i$ be the image of ${\mathbf b}_i$ ($i=1\dots n$) under the homomorphism $V\rightarrow V/W$.
From Lemma \ref{mats} it follows that there are two disjoint index sets $I_1,I_2\subset \{1\dots n\}$ such that
$\{{\mathbf p}_i|i\in I_1\}$ and $\{{\mathbf p}_i|i\in I_2\}$ are bases of the same subspace $T\leq W$.
Let $j$ be an element in $I_1$, and let ${\mathbf x}=\sum_{i=1}^{n} \lambda_i{\mathbf b}_i$
be an element in $H_j\leq V$ which is not covered by $H_i+{\mathbf v}_i$ for $i\neq j,~1\leq i\leq k$.
Since the hyperplanes $H_l$~$(1\leq l\leq n)$ does not cover
${\mathbf x}$ for $l\neq j$ it follows that $\lambda_l\neq 0$ for
all $l\neq j$.
Let ${\mathbf y}=\sum_{i=1}^n \lambda_i{\mathbf p}_i$ and ${\mathbf y}_1=\sum_{i\in I_1\cup I_2} \lambda_i{\mathbf p}_i$.
Lemma \ref{brumm} implies that ${\mathbf y}_1$ can be written as a nowhere zero linear combination of the vectors
${\mathbf p}_i$~$(i\in I_1\cup I_2)$ and thus ${\mathbf y}$ can be written in the form
$\sum_{i=1}^{n} \mu_i{\mathbf p}_i$ where $\mu_i\neq 0$ for $1\leq i\leq n$.
Let ${\mathbf z}= \sum_{i=1}^{n} \mu_i{\mathbf b}_i$. The vector ${\mathbf z}$ is a preimage of ${\mathbf y}$
under the homomorphism $V\rightarrow W$ and so ${\mathbf z}-{\mathbf v}\in W$. Since ${\mathbf z}$ is a
nowhere zero vector in the basis $B$ we have that it is not contained in $H_1,H_2,\dots,H_n$. Let
$t>n$ be a number for which $H_t+{\mathbf v}_t$ contains ${\mathbf z}$. By definition of $W$, $H_t+{\mathbf v}_t$ contains
the set ${\mathbf z}+W$. This contradicts the assumption that ${\mathbf v}$ is covered only by $H_j$.
\end{proof}

Note that the condition on $q$ was hidden in Lemma \ref{brumm} when we used a result of \cite{AT}.
This means that the choosability version of the Alon-Jaeger-Tarsi conjecture would imply the analogue statement for an
arbitrary prime number bigger than $3$. It can also be seen (from the previous proof)
that the following weak conjecture implies
$\varepsilon_p>0$ if $p>2$.

\medskip
\noindent{\bf Weak conjecture}
{\it For every prime $p>2$ there is a constant $c_2(p)$ depending only on $p$ such that if $B_1,B_2,\dots ,B_{c_2(p)}$
are bases of the $GF(p)$ vectorspace $V$ then all elements of $V$ can be written as a nowhere zero linear combination
of the elements of the union (as multisets) of the previous basises.}
\medskip

The next result shows that the weak conjecture is equivalent with $\varepsilon_p>0$.

\begin{lem}
If $\varepsilon_p>0$ then the weak conjecture holds for $p$ with any
$c_2(p)=k>{{1+\varepsilon_p}\over{\varepsilon_p}}$ .
\end{lem}

\begin{proof}
We go by contradiction. Assume that the weak conjecture is
not true with $c_2(p)=k$. Let $n$ be the minimal dimension where the conjecture is false (with $c(p)=k$) and assume
that the bases $B_1,B_2,\dots,B_k$ are forming a counter example in the $n$ dimensional space $V$.
Let $M$ be an $n$ by $nk$ matrix whose columns are
the vectors from the previous bases. According to our assumption there is a vector ${\mathbf v}$ such that
there is no nowhere zero vector ${\mathbf x}$
with $M{\mathbf x}={\mathbf v}$.
Let say that an index set $I\subseteq \{1...nk\}$ is a blocking set if for all ${\mathbf x}\in {\rm GF}(p)^{nk}$ with
$M{\mathbf x}={\mathbf v}$
there is a $j\in I$ such that the $j$-th coordinate of ${\mathbf x}$ is zero.
Let $I$ be a minimal blocking set. First we prove by contradiction that $I=\{1...nk\}$.
Assume that $P=\{1..nk\}\setminus I$ is not empty. Let $j$ be an element of $P$, let ${\mathbf y}$ be the $j$-th
column of $M$ and let $W$ be the factor space $V/\langle{\mathbf y}\rangle$. Let $P_1,P_2,\dots,P_k$ be the
images of the bases $B_1,B_2,\dots,B_k$ under the homomorphism $V\rightarrow W$. It is clear that each
$P_i$ contains a basis for $W$ and by minimality of $n$ it follows that each vector ${\mathbf x}\in W$ is a
nowhere zero linear combination of the elements in $P_1,P_2,\dots,P_k$. In particular the image of ${\mathbf v}$ can
be written as such a nowhere zero combination. It means that there is a vector ${\mathbf x}\in {\rm GF}(p)^{nk}$
for which $M{\mathbf x}={\mathbf v}$ and all but the $j$-th coordinate of ${\mathbf x}$ are not zero.
It contradicts the assumption that $I$ is a blocking set.
Now we have that $\{1..nk\}$ is a minimal blocking set and thus for each $j\in \{1..nk\}$ there is a vector
${\mathbf x}_j\in {\rm GF}(p)^{nk}$ such that all but the $j$-th coordinate of ${\mathbf x}_j$ are not zero and
$M{\mathbf x}_j={\mathbf v}$. Let $U$ be the affine hyperplane consisting of all ${\mathbf x}$ for which
$M{\mathbf x}={\mathbf v}$.
For all $j\in \{1..nk\}$ let $H_j\leq U$ be the affine hyperplane consisting of those
elements ${\mathbf x}$ whose $j$-th coordinate is zero.
Now the affine space $U$ is covered irredundantly by the affine hyperplanes $H_j$.
Since ${\rm dim}(U)=n(k-1)$ it follows that ${k\over{k-1}}\geq 1+\varepsilon_p$.
\end{proof}

\bigskip

\section{colorings and flows}

In this section we outline the relation between colorings, flows and hyperplane coverings.
Let $G$ be a finite, loopless graph with vertex set $V(G)$ and edge set $E(G)$.
Let $q$ be a prime power and let $W$ be the vector space of all functions $V(G)\rightarrow {\rm GF}(q)$.
For two functions $f,g\in W$ we define their scalar product by
$$(f,g)=\sum_{v\in V(G)}f(v)g(v)$$
We associate a vector ${\mathbf v}_e\in W$ to each edge $e\in E(G)$ such that
${\mathbf v}_e$ takes $1$ and $-1$ at the two different
endpoints of $e$, and it takes $0$ everywhere else.

\begin{lem}
$G$ is colorable with $q$ colors if and only if the orthogonal spaces of the vectors ${\mathbf v}_e$ do not
cover the whole space $W$.
\end{lem}

\begin{proof}
We can think of $W$ as the set of all possible (not necessary proper) colorings of $G$.
It is clear that a vector ${\mathbf v}\in W$ is orthogonal to ${\mathbf v}_e$ for some $e\in E(G)$ if and only
if ${\mathbf v}$ takes the same value at the endpoints of $e$. It means that $G$ has a proper coloring with $q$ colors
if and only if there is a vector ${\mathbf v}\in W$ which is not contained in any of the spaces
${\mathbf v}_e^{\bot}$.
\end{proof}

Combining the previous lemma with Theorem \ref{chariz} one gets the following peculiar
characterization of colorability.

\medskip

\begin{proposition}
If $q$ is an odd prime then $G$ can be colored by $q$ colors if
and only if there is a vector ${\mathbf v}\in W$ such that the
number of zero-one combinations of the vectors ${\mathbf v}_e$
resulting ${\mathbf v}$ is odd.
\end{proposition}

\medskip

Note that the space ${\mathbf v}_e^{\bot}$ depends only on the one dimensional space spanned by ${\mathbf v}_e$.
It means that the vectors ${\mathbf v}_e$ can be replaced by any nonzero representative from their one dimensional
spaces, which gives an even stronger version of the previous proposition. We also note that the "if" direction
remains true if we delete the condition that $q$ is a prime number.

Let $G=(V,E)$ be a directed graph and let $A$ be an Abelian group. An $A$-flow on $G$ is a function $f:E\rightarrow A$
such that for all $v\in V$
\begin{equation*}
\sum_{e\in \delta^+ (v)}f(e)=\sum_{e\in \delta^- (v)}f(e).
\end{equation*}
where $\delta^+ (v)$ denote the set out going edges and $\delta^- (v)$ denote the set of in coming edges.
If $f(e)\neq 0$ for all $e\in E$ then $f$ is called a nowhere zero flow.
Clearly the existence of a nowhere zero flow on $G$ is independent of the orientation of $G$.
If $G$ is undirected we will say that it admits a nowhere zero $A$-flow if some (an thus every) orientation of it
admits a nowhere zero $A$-flow.
Let $G$ be a fixed graph with a fixed direction and consider the set $B$ of all possible flows on $G$. It is
clear that $B$ is a subgroup of the direct product $A^E$ and one can prove easily that
$B\simeq A^{|E|-|V|+m}$ where $m$ denotes the number of connected components of $G$. For each edge $e$ there is a subgroup
$B_e\leq B$ consisting of those flows which vanish on $e$. Clearly, $G$ has a nowhere zero flow if and only if the
subgroups $B_e$ do not cover the group $B$. Moreover, it is also clear that the intersection of the
subgroups $B_e$ is trivial. It means in particular that if $G$ is a graph which is "edge-minimal" respect to the
property of having no nowhere zero flow (i.e. $G$ has no nowhere zero flow, but after deleting any edge, the resulting
graph always has one) then the number of edges is less than $g_1(|B|)$ where $g_1$ is the function defined in the introduction.
Note that if $A$ has a finite field structure, then the group $B$ can be regarded as
a vectorspace over $A$ with hyperplane system $\{B_e|e\in E\}$.

\section{hierarchy of conjectures}

\begin{picture}(100,100)(10,10)
\put (0,80){case $p>3$}
\put (0,50){$c_2=2$}
\put (30,53){\vector(1,0){20}}
\put (55,50){$\varepsilon_p\geq {\rm log}_2(p)-1$}
\put (130,53){\vector(1,0){20}}
\put (155,50){$\varepsilon_p>1$}
\put (185,53){\vector(1,0){20}}
\put (210,50){C-AJT}
\put (243,53){\vector(1,0){20}}
\put (268,50){$\varepsilon_p\geq 0.5$}
\put (280,48){\vector(0,-1){20}}
\put (268,20){$\varepsilon_p>0$}
\put (264,25){\vector(-1,0){20}}
\put (244,21){\vector(1,0){20}}
\put (230,20){W}
\put (208,23){\vector(1,0){20}}
\put (190,20){AB}
\put (223,60){\vector(1,1){20}}
\put (240,85){AJT}
\end{picture}

\begin{picture}(50,50)(10,10)
\put(0,40){case $p=3$}
\put(0,0){$c_2=2$}
\put(30,3){\vector(1,0){20}}
\put(55,0){$\varepsilon_3\geq {\rm log}_2(3)$}
\put(110,3){\vector(1,0){20}}
\put(133,0){$\varepsilon_3>0$}
\put(163,1){\vector(1,0){20}}
\put(183,5){\vector(-1,0){20}}
\put(187,0){W}
\put(199,1){\vector(1,0){20}}
\put(219,5){\vector(-1,0){20}}
\put(222,0){AB}
\put(240,3){\vector(1,0){20}}
\put(262,0){WT}
\end{picture}

\begin{picture}(10,10)(10,10)
\end{picture}

\bigskip

\medskip
\begin{tabular}{ll}
\begin{tabular}{l}
AJT\\
C-AJT\\
AB\\
W\\
WT\\
\end{tabular}
&
\begin{tabular}{l}
Alon-Jaeger-Tarsi conjecture\\
choosability version of AJT\\
additive basis conjecture\\
weak conjecture\\
weak three flow conjecture\\
\end{tabular}
\end{tabular}

\bigskip

{\bf Aknowledgements}
The author thanks N. Alon, P.P. P\'alfy, L. Pyber and C. Szegedy for their kind help and helpful remarks.

\end{document}